\documentclass{article}[10]

\usepackage[a4paper, total={6in, 8in}]{geometry}
\usepackage{amsmath}
\usepackage{amsmath, amsthm}
\usepackage{amsopn}
\usepackage{graphicx, amssymb, array, amscd, amsfonts, color}
\usepackage{latexsym}

\theoremstyle{definition}

\theoremstyle{remark}

\theoremstyle{problem}
\newtheorem{problem}{Problem}
\date{}

\title{Chebyshev Differential Quadrature for Numerical Solutions of Higher Order Singular Perturbation Problems}
\author{ G\"ulsemay   Yi\u git$^{1}$ and Mustafa Bayram$^{2}$\\ 
  \small $^{1}$Department of Basic Sciences, Istanbul Kemerburgaz University, Istanbul, Turkey \\
		\small $^{2}$Department of Computer Engineering, Istanbul Geli\c sim University, Istanbul, Turkey \\
  \small $^{1}$gulsemay.yigit@kemerburgaz.edu.tr \\
		\small $^{2}$mbayram@gelisim.edu.tr
}

\begin{document}
\maketitle

{\abstract 
In this study linear and nonlinear higher order singularly perturbed problems are examined by a numerical approach, the differential quadrature method. Here, the main idea is using Chebyshev polynomials to acquire the weighting coefficient matrix which is necessary to get numerical results. Following this, different class of perturbation problems are considered as test problems. Then, all results are shown in tables and also comparison between numerical and exact solution shows the accuracy and effectiveness of the presented method.}  
\\
\begin{flushleft}
\\ 
\end{flushleft}
\section{Introduction}
The analytical and numerical methods to solve singular perturbation problems have been widely used in many fields of fluid dynamics, reaction-diffusion processes, particle physics, and combustion processes. These types of problems are represented by differential equations including $\epsilon$  which is assumed to be a small parameter and solutions of the problems have non-uniform behaviour when the parameter $\epsilon\rightarrow 0$. For detailed information one can consider books given in [1-3]. \\
\\
In this study, it is considered singularly perturbed various kind of higher order problems.  These problems are classified by setting   $\epsilon=0$  . If the order of the higher order problem is reduced by one, the problem becomes convection-diffusion type and if the order is reduced by two, it becomes reaction-diffusion type [4]. In this context, it is developed an efficient numerical method to solve both linear and nonlinear  higher order problems. Then, to show the efficiency of the method, different types of problems have been solved and the solutions are compared with analytical or asymptotic methods depending on the small parameter $\epsilon$ . \\
\\
Many studies have been carried out on singular perturbation including both theoretical and numerical. Howes [5] has considered third order singularly perturbation problems and investigate the existence of problems. The existence of solution of third order problem is also examined by Yao and Feng [6]. Other than perturbation methods, many numerical techniques for solving such kind problems available in literature. Boundary value technique is one of the popular methods for solving such problems. It is the method for obtainig numerical results for higher order singularly perturbed differential equations by transforming the problem to a system of a first differential equation and a second order singularly perturbed equation. [7-8] Such problems with discontinuous source term  are also considered by asymptotic finite element method [9]. Cui and Geng, studied on third order singularly perturbed problem using an analytical method with asymptotic expansion to have solution for the problem [10]. Kumar and Tiwari considered reaction-diffusion kind boundary value problems by using initial-value technique. According to the idea, the  third order perturbation problem is transformed into three unperturbed initial value problems to have solutions numerically [4]. Phaneendra, Reddy and Soujanya solved such kind of problems by using Fitted Numerov Method by transforming third order boundary value problem to the second order. [11]. \\
\\
Differential Quadrature is a numerical approximation technique to solve ordinary or partial 
differential equations and presented by Bellman and his associates [12-13]. The goal of the technique is to present a new perspective for previous numerical techniques in solving problems. Since then, this method has been used in a wide range of applications. \\
\\
The DQ method is similar to the integral quadrature. According to this method, derivative of a function is defined as a weighted linear sum of the  function values at all grid points related to the used direction. So, there occurs the term weighting coefficients and to obtain these coefficients, generally polynomials are chosen as test functions which can be obtained by polynomial approximation theory.\\
\\Firstly, Bellman proposed two polynomial-based methods to compute the weighting coefficients of the first order derivative. In the first method, power function is used as test function. And the second one depends upon choosing test function as Legendre polynomials [12-13]. Soon afterwards, many different approaches have been presented and applied to many different kind of engineering problems such as using Lagrange polynomials, Spline functions, Fourier expansion.  Quan and Chang [14-15] used Legendre interpolated polynomials as a test functions then obtained explicit formulation to find the weighting coefficients. In this study, we used Chebyshev polynomial approximations to obtain numeric solutions. When the method is applied for the derivatives, differential equation is reduced to linear system of equation, with the implementation of boundary or initial conditions, matrix equation can be solved to obtain the desired solution.\\
\\
One of the most popular improvements related to method made by Shu [16]. He generated a useful way combining Bellman’s and Quang-Chang’s approaches. In Shu’s book, first and higher order derivative formulations are analyzed in detail based on polynomial approaches using different kind of grid points. 
\\
Many applications on this method, available in literature such as fluid mechanics, structural mechanical analysis, chemical processes, biosciences,  etc. And throughout these areas, it can be seen that this technique gives accurately solutions with time saving computations [17]. Mingle [18-19], applied this method to nonlinear diffusion equation first and he gives an approach for implementing boundary conditions for the differential equations. Civan and Spliepcevich applied this method to both Poisson equation [20] and to multi-dimensional problems [21]. 
\\
 Y\"ucel and Sar{\i}, applied the DQ method for solving singular boundary value problems. It is shown that numerical examples for some specific boundary conditions, method gives highly accurate solutions [22]. Also, Y\"ucel considered Sturm-Liouville problem by using Lagrange polynomials and and Fourier series expansion to evaluate the weighting coefficients, results were compared with the previous works and it was shown that the method gives more accurate results [23]. Saka et al. considered equal width equation (EW) by using three methods including cosine expansion based differential quadrature.  It is observed that numerical solutions of the given problem can be obtained with high accuracy by using the differential quadrature method [24]. Sar{\i} applied the quadrature method for solving two point singular perturbation analyzing linear and nonlinear second order boundary value problems and solutions were compared with the previous works [25]. Alper et al, considered many engineering problems using both Spline functions and polynomials as test functions. For time discretization they used Runge-Kutta methods and also stability analysis for the method is investigated [26-28].
\\

\section{The Differential Quadrature Method (DQM) }
In this section differential Quadrature Method is analyzed by using polynomials as test functions and to get weighting coefficients Chebyshev polynomials are used. \\
\\ 
Consider a sufficiently smooth function      $f(x)$      on a closed interval     $[a,b]$        Derivative of the function  at a grid point   $x_i$   , is approximated by a linear sum of all functional values in the whole domain and the quadrature formula for first derivative is given as follows:

\begin{equation} \label{eq:1}
f_x(x_i)={\dfrac{df}{dx}}{\biggr{|}_{x_i}} = \sum_{j=1}^{N} a{_{ij}} f(x_j)~~~~i=1,2,\cdots,N
\end{equation}
\\
where   $a{_{ij}}$     represents the weighting coefficients, $N$ is the number of grid points. Here, the main idea is to determine weighting coefficients. [16]
\\
To generalize the idea, the n-th order derivative is given as,  
\\
\begin{equation} \label{eq:2}
 f_x^{(n)}(x_i)={\dfrac{df}{dx}}{\biggr{|}_{x_i}} = \sum_{j=1}^{N} w{_{ij}^{(n)}} f(x_j)~~~~i=1,2,\cdots,N
\end{equation}
where     $w{_{ij}^{(n)}}$   represents the weighting coefficients, N is the number of grid points. Here, the main idea is to determine weighting coefficients [16].
\\
Here, we choose grid points as the Chebyshev collation points defined as
\begin{equation} \label{eq:3}
 x_i=\cos(\theta_i),~~~~ \theta=\dfrac{i\pi}{N}, ~~~~i=1,2,\cdots,N
\end{equation}
\\
which is applicable for only interval [1,-1]. If the problem is given on interval [a,b] to obtain $x_i$     following transformation is used [16].
\begin{equation} \label{eq:4}
 x_i=\dfrac{b-a}{2} (1-\xi_i)+a.
\end{equation}
If Lagrange interpolating polynomials are used, 
\begin{equation} \label{eq:5}
 r_k (x)={\dfrac{K(x)}{(x-x_k) K^\prime (x_k)}}, ~~~~i=1,2,\cdots,N
\end{equation}
where $r_k (x)$ represents the test function and 
\begin{equation} \label{eq:6}
 K(x)=(x-x_1)(x-x_1)(x-x_2)\cdots(x-x_N),~~~~ K^\prime (x)=\prod\limits_{k=1,k\neq j}^{N}(x_j-x_k).
\end{equation}
Then, with some linear algebra, the weighting coefficients for first order derivative given in (\ref{eq:1})
become as follows [16]:

\begin{equation} \label{eq:7}
 a_{ij} (x)={\dfrac{K^\prime (x_i)}{(x_i-x_j) K^\prime (x_j)}},~~~j \neq i.
\end{equation}\\
For the entries on main diagonal, the following relation becomes,
\begin{equation} \label{eq:8}
 a_{ii} (x)=\dfrac{{K''} (x_i)}{(x_i-x_j) 2K^\prime (x_i)}.
\end{equation} \\
Since, a linear vector space can be spanned by different kind of bases, for diagonal entries following formula can be used which is obtained by using power functions, $x^{k-1},~~~k=1,2,\cdots,N $ when $k=1$,

\begin{equation} \label{eq:9}
\sum\limits_{j=1}^{N}a{_{ij}}=0,~~~a_{ii}=-\sum\limits_{{j=1},j \neq i}^{N}a{_{ij}}.
\end{equation}\\
Here, to obtain quadrature solutions of the problems it is used Chebyshev polynomial as a basis. The Chebyshev polynomial of the first kind $T_k(x)=\cos(k\theta),~\theta=\arccos(x)$ and Chebyshev approximation polynomial is given as  
\begin{equation} \label{eq:10}
f(x)=\sum\limits_{j=1}^{N}c_k T_k(x).
\end{equation}
where $c_k$ is the Chebyshev approximation coefficient. Then, for each grid point $n-$th order derivative approximation is obtained as,
\begin{equation} \label{eq:11}
 f_x^{(n)}(x_i)= \sum\limits_{j=1}^{N} w{_{ij}^{(n)}} f(x_j)~~~~i=1,2,\cdots,N.
\end{equation}
where $w{_{ij}^{(n)}}$ represents the $n-$th order weighting coefficient matrix, and there are suitable ways in literature, considered to find first derivative matrix, as follows [16],[29]:
\begin{eqnarray*} \label{eq:12}
&& w_{ij} =\dfrac{{{\overline{c_i}}}(-1)^{i+j}}{{{\overline{c_j}}}(x_i-x_j) },~~~0\leq i,j\leq N, i \neq j,\nonumber\\
&& w_{ii} =-\dfrac{x_i}{2(1-x^{2}_{i}) },~~~1\leq i\leq N-1, \nonumber\\
&& w_{00} =-w_{NN}={\dfrac{2N^2+1}{6} }.
\end{eqnarray*}
where $\overline{c_{0}}={\overline{c_{N}}}=2$ and ${\overline{c_j}}=1,~~~0\leq i\leq N-1.$
\\
\section{Applications to Higher Order Singular Perturbed Problems}
In this section,it is investigated different type of higher order perturbation problems, to obtain Chebyshev quadrature solutions.
\\ 
\subsection{Third Order Reaction Diffusion Problems}\label{sec:3.1}
At first, we consider reaction diffusion type BVP. Such problems are given by third order ordinary differential equation with a small parameter [4],
\begin{equation} \label{eq:13}
\epsilon y^{\prime\prime\prime}(x)+p(x)y^{\prime}(x)+q(x)y(x)=f(x),~~~x\in [a,b],
\end{equation}
with subject to the following boundary conditions
 \begin{equation*} 
y(a)=a_0,~~~y(b)=b_0 ~~~y^{\prime}(a)=c_0.
\end{equation*}
where $\epsilon>0$ is a small parameter and $a_0,b_0,c_0$ are known constant values and $p(x),q(x)$ and $f(x)$ are differentiable functions on $a\leq x\leq b$. To obtain the approximated solution, we apply the method for each grid points, as follows,
\begin{equation} \label{eq:14}
\epsilon y^{\prime\prime\prime}(x_i)+p(x_i)y^{\prime}(x_i)+q(x_i)y(x_i)=f(x_i).
\end{equation}
Then, derivatives are replaced by the Differential Quadrature equality,
\begin{equation} \label{eq:15}
\epsilon \sum\limits_{j=1}^{N} w{_{ij}^{(3)}} y{(x_j)}+p(x_i)\sum\limits_{j=1}^{N} w{_{ij}^{(1)}}y(x_j)+q(x_i)y(x_i)=f(x_i),
\end{equation}
or by simplifying
\begin{equation} \label{eq:16}
 \sum\limits_{j=1}^{N} \left\{ \epsilon w{_{ij}^{(3)}}+p(x_i) w{_{ij}^{(1)}}+q(x_i)  \right\}y{(x_j)}=f(x_i).
\end{equation}
where $w{_{ij}^{(3)}}$ is the third order and $w{_{ij}^{(1)}}$ represents the first order weighting coefficients which is obtained by using formulations given by equation \eqref{eq:12}.Then, the equation is represented by system of eqautions given by,
\begin{equation} \label{eq:17}
 \left[\epsilon w{_{ij}^{(3)}}+p(x_i) w{_{ij}^{(1)}}+q(x_i)I \right]  [y{(x_j)}]=[f(x_i)].
\end{equation}
For simplicity, the equation \eqref{eq:17} is represented by
\begin{equation} \label{eq:18}
 \left[A_{ij}\right] \left[y(x_j)\right]=\left[f(x_i)\right].
\end{equation}\\
The method is also applied for the boundary conditions as follows,
\begin{eqnarray} \label{eq:19}
&& y(a)=y(x_1)=a_0, ~y(b)=y(x_N)=b_0, \nonumber\\ 
&& {y^{\prime} (a)=y^\prime (x_1) =\sum\limits_{j=1}^{N} w{_{1j}^{(1)}} y(x_j)}=c_0. 
\end{eqnarray}
Here, the matrix, $[A_{ij}]$  is the weighting coefficient matrix, including third and first order derivatives and the type of this matrix is skew centrosymmetric. That is,determinant of the matrix is zero and inverse does not exist. To obtain an accurate solution, we should replace any rows with boundary conditions. So the matrix equation becomes,
\begin{equation*} 
\begin{bmatrix}
    1 & 0 & 0 &\cdots& 0 \\
    w{_{11}^{(1)}} & w{_{12}^{(1)}}& w{_{13}^{(1)}}&\cdots&w{_{1N}^{(1)}} \\
		A{_{31}} & A{_{32}}& A{_{33}}&\cdots& A{_{3N}} \\
		\vdots & & & \ddots&\\
		A{_{(N-1)1}} & A{_{(N-1)2}}& A{_{(N-1)3}}&\cdots& A{_{(N-1)N}} \\
		0 & 0 & 0 &\cdots& 1 \\
  \end{bmatrix}
	\begin{bmatrix}
    y_1  \\
    y_2 \\
		y_3 \\
    \vdots\\	
		y_{N-1} \\
		y_N \\
  \end{bmatrix}
	=
	\begin{bmatrix}
    a_0 \\
    c_0 \\
		f_3 \\
    \vdots\\	
		f_{N-1} \\
		b_0 \\
  \end{bmatrix}.
	\end{equation*}	
After implementing the boundary conditions, the system is computed to obtain numerical solutions  $y_i$. Here, it is used LU decomposition method.

\subsection{Fourth Order Convection Diffusion Problems}\label{sec:3.2}
In this case, it will be considered convection diffusion type singular perturbed boundary value problem. Such problems are given by ordinary differential equation as follows:  
\begin{equation} \label{eq:20}
-\epsilon y^{\prime\prime\prime\prime}(x)-p(x)y^{\prime\prime\prime} +q(x)y^{\prime\prime}(x)-r(x)y(x)=-f(x),~~~x\in (0,1),
\end{equation}
with subject to the following boundary conditons,
 \begin{equation*} 
y(0)=a_0,~~~y(1)=b_0, ~~~y^{\prime\prime}(0)=-c_0,~~~y^{\prime\prime}(1)=-d_0.
\end{equation*}
where $\epsilon>0$ is a small parameter and $p(x), q(x), r(x)$ and $f(x)$ are differentiable functions on $(0,1)$. Using the same idea, DQ discretization for fourth order type problem is of the form
\begin{equation} \label{eq:21}
-\epsilon \sum\limits_{j=1}^{N} w{_{ij}^{(4)}} y{(x_j)}-p(x_i)\sum\limits_{j=1}^{N} w{_{ij}^{(3)}} y{(x_j)}q(x_i)\sum\limits_{j=1}^{N} w{_{ij}^{(2)}}y(x_j)-r(x_i)y(x_i)=-f(x_i),
\end{equation}
Then, boundary conditions are implemented to the system. After all, we obtain  $N\times N$ system of  equations, and given by
\begin{equation} \label{eq:22}
 [A_{ij}] [y(x_j)]=[f(x_i)].
\end{equation}
Here, the matrix,   $  [A_{ij}]$  s the weighting coefficient matrix, including fourth,third and second order derivatives. Finally, the system can be solved by suitable matrix solver to obtain numerical solutions   $y_i$. 
\\
\section{Numerical Illustrations}
To show the accuracy of the method, approximated error is computed by using  $L_2$     error norm which is given by,
\begin{equation} \label{eq:23}
 L_2={|U_{ex}-U_{N}|}={{\left[\sum\limits_{j=1}^{N} {\left|{(U_{ex})}_j-{(U_{N})}_j\right|}^2\right]}}^{1/2}  i=1,2,\cdots,N
\end{equation}
and maximum error norm given by
\begin{equation} \label{eq:24}
 L_\infty={|U_{ex}-U_{N}|}_{\infty}=\underset{j}{max}{\left|{(U_{ex})}_j-{(U_{N})}_j\right|}~~~   i=1,2,\cdots,N.
\end{equation}
where $U_{ex}$ and $U_{N}$ represents the analytical and numerical solutions. We have used differential quadrature derivative approximations to obtain numerical results. Now, five test problems are solved to illustrate the efficiency of the presented method and all the results in terms of indicated norms are given in tables.  
\\
\newpage
\begin{problem}
Reaction Diffusion Type 

\end{problem}

\begin{equation} \label{eq:25}
\epsilon y^{\prime\prime\prime}+4y^{\prime}-4y=x^2.
\end{equation}
 \begin{equation*} 
y(0)=0.5,~~~y(1)=1.47 ~~~y^{\prime}(0)=0.5
\end{equation*}
The problem in \eqref{eq:25} is examined by initial value technique in [4].
\begin{center}
{\textbf{Table 1.~~}Error Norms for Problem 1 }
	\label{tab:Error Values for Problem 1}
    \begin{tabular}{  l  l  l  l  l   }
    \hline
		\\
    \textbf{Error Norms} & \textbf{$~~N~~$} & $~~~~\epsilon=10^{-1}~~$  & $~~~~\epsilon=10^{-2}~~$ & $~~~~\epsilon=10^{-3}$~~  \\  \hline
    
		\\
		$$ &~~10~~ & ~~0.3395E-04 & ~~0.9940E-02& ~~0.2397E-04      \\ 
    $L_2$ & ~~20 ~~& ~~0.6561E-13 & ~~0.3351E-04& ~~0.4381E-04    \\ 
    $$ & ~~50~~& ~~0.6224E-13 & ~~0.2081E-13& ~~0.6962E-07  \\

		\hline
		\\
		$$ &~~10~~ & ~~0.2024E-04 & ~~0.5844E-04& ~~0.1576E-04      \\ 
    $L_\infty$ & ~~20 ~~& ~~0.2607E-13 & ~~0.4510E-04& ~~0.2019E-04    \\ 
    $$ & ~~50~~& ~~0.1629E-13 & ~~0.5424E-14& ~~0.6224E-06  \\ 
		
			\hline

		\\
    \end{tabular}

\end{center}
Table 1 shows the error for increasing number of grid points between exact solution and numerical solutions. It can be seen that, when using 10 grid points and $\epsilon \longrightarrow 0$ numerical solutions are less stable when compared with case for using grid points 50.

\begin{problem}
Reaction Diffusion Type 
\end{problem}

\begin{equation} \label{eq:26}
\epsilon y^{\prime\prime\prime}+\left(1+\dfrac{x}{2}\right)y^{\prime}-\dfrac{1}{2}y=0.
\end{equation}

\begin{equation*} 
y(0)=0.6,~~~y(1)=0.9 ~~~y^{\prime}(0)=0.23
\end{equation*}\\
The problem in \eqref{eq:26} is examined by initial value technique in [4]
\\               
\begin{center}

{\textbf{Table 2.~~}Error Norms for Problem 2 }
	\label{tab:Error Values for Problem 2}
    \begin{tabular}{  l  l  l  l  l   }
    \hline
		\\
    \textbf{Error Norms} & \textbf{$~~N~~$} & $~~~~\epsilon=10^{-1}~~$  & $~~~~\epsilon=10^{-2}~~$ & $~~~~\epsilon=10^{-3}$~~  \\  \hline
    
		\\
		$$ &~~10~~ & ~~0.1923E-05 & ~~0.7964E-02& ~~0.3463E-02      \\ 
    $L_2$ & ~~20 ~~& ~~0.5529E-15 & ~~0.1407E-08& ~~0.3660E-02    \\ 
    $$ & ~~50~~& ~~0.3388E-15& ~~0.1876E-13& ~~0.6962E-12  \\

		\hline
		\\
		$$ &~~10~~ & ~~0.1095E-05 & ~~0.4566E-02& ~~0.2198E-02      \\ 
    $L_\infty$ & ~~20 ~~& ~~0.2132E-15 & ~~0.3809E-04& ~~0.2986E-02    \\ 
    $$ & ~~50~~& ~~0.9173E-15 & ~~0.3748E-14& ~~0.2406E-13  \\ 
		
			\hline

		\\
    \end{tabular}

\end{center}
Table 2 shows the error for increasing number of grid points. It can be seen that, when $\epsilon \longrightarrow 0$ numerical solutions are less stable.
\\
                
\begin{problem}
Fourth Order Convection Diffusion Type 
\end{problem}

\begin{equation} \label{eq:27}
-\epsilon y^{\prime\prime\prime\prime}-4y^{\prime\prime\prime}=1.~~~~0<x<1,
\end{equation}

\begin{equation*} 
y(0)=1,~~~y(1)=1 ~~~y^{\prime\prime}(0)=-1,~~~y^{\prime\prime}(1)=-1.
\end{equation*}\\
The example given by \eqref{eq:27} is solved by Boundary Value Technique in [30].\\
\begin{center}
{\textbf{Table 3.~~}Error Norms for Problem 3 }
	\label{tab:Error Values for Problem 3}
    \begin{tabular}{  l  l  l  l  l   }
    \hline
		\\
    \textbf{Error Norms} & \textbf{$~~N~~$} & $~~~~\epsilon=10^{-1}~~$  & $~~~~\epsilon=10^{-2}~~$ & $~~~~\epsilon=10^{-3}$~~  \\  \hline
    
		\\
		$$ &~~10~~ & ~~0.8992E-03 & ~~0.9301E-03& ~~0.9332E-03      \\ 
    $L_2$ & ~~20 ~~& ~~0.1311E-04 & ~~0.3159E-04& ~~0.7893E-03    \\ 
    $$ & ~~50~~& ~~0.5141E-09 & ~~0.1361E-05& ~~0.3116E-03  \\

		\hline
		\\
		$$ &~~10~~ & ~~0.5275E-03 & ~~0.5497E-03& ~~0.5402E-03      \\ 
    $L_\infty$ & ~~20 ~~& ~~0.5221E-05 & ~~0.8033E-04& ~~0.3102E-04    \\ 
    $$ & ~~50~~& ~~0.2580E-09 & ~~0.3222E-06& ~~0.7445E-06  \\ 
		
			\hline

		\\
    \end{tabular}

\end{center}
Table 3 shows the error for increasing number of grid points. As the grid points increase more accurate solutions is obtained. But, when $\epsilon \longrightarrow 0$ numerical solutions are approximately same. 
\\           

\begin{problem}
Nonlinear Boundary Value Problem
\end{problem}

\begin{equation} \label{eq:28}
\epsilon y^{\prime\prime\prime}+y^{\prime\prime}+\epsilon y^{\prime}(y^{\prime}+2)=1.~~~~0\leq<x\leq<\dfrac{\pi}{2},
\end{equation}

\begin{equation*} 
y(0)=0,~~~y(\dfrac{\pi}{2})=1-\dfrac{\epsilon}{3} ~~~y^{\prime}(\dfrac{\pi}{2})=-1+\dfrac{\epsilon}{4}.
\end{equation*}\\
The problem given by \eqref{eq:28} is analyzed asymptotically in [1].
\begin{center}
{\textbf{Table 4.~~}Error Norms for Problem 4 }
	\label{tab:Error Values for Problem 4}
    \begin{tabular}{  l  l  l  l  l   }
    \hline
		\\
    \textbf{Error Norms} & \textbf{$~~N~~$} & $~~~~\epsilon=10^{-1}~~$  & $~~~~\epsilon=10^{-2}~~$ & $~~~~\epsilon=10^{-3}$~~ \\  \hline
    
		\\
		$$ &~~10~~ & ~~0.8992E-03 & ~~0.9301E-03& ~~0.9332E-03      \\ 
    $L_2$ & ~~20 ~~& ~~0.1311E-04 & ~~0.3159E-04& ~~0.7893E-03    \\ 
    $$ & ~~50~~& ~~0.5141E-09 & ~~0.1361E-05& ~~0.3116E-03  \\

		\hline
		\\
		$$ &~~10~~ & ~~0.5275E-03 & ~~0.5497E-03& ~~0.5402E-03      \\ 
    $L_\infty$ & ~~20 ~~& ~~0.5221E-05 & ~~0.8033E-04& ~~0.3102E-04    \\ 
    $$ & ~~50~~& ~~0.2580E-09 & ~~0.3222E-06& ~~0.7445E-06  \\ 
		
			\hline

		\\
    \end{tabular}

\end{center}
Table 4 shows the error for increasing number of grid points.Because of the asymptotic behaviour of the problem it is observed that for all grid points results have approximately same accuracy. 
\\
\begin{problem}
Nonlinear Initial Value Problem
\end{problem}   
\begin{equation} \label{eq:29}
\epsilon y^{\prime\prime\prime}+y^{\prime\prime}+\epsilon \left[{(y^{\prime})}^2+y\right]=\epsilon e^{-2x},x\geq0,
\end{equation}

\begin{equation*} 
y(0)=2,~~~y^{\prime}(0)=-1~~~y^{\prime}(0)=1.
\end{equation*}
The problem given by \eqref{eq:29} is analyzed asymptotically in [1].

\begin{center}
{\textbf{Table 5.~~}Error Norms for Problem 5 }
	\label{tab:Error Values for Problem 5}
    \begin{tabular}{  l  l  l  l  l   }
    \hline
		\\
    \textbf{Error Norms} & \textbf{$~~N~~$} & $~~~~\epsilon=10^{-1}~~$  & $~~~~\epsilon=10^{-2}~~$ & $~~~~\epsilon=10^{-3}$~~  \\  \hline
    
		\\
		$$ &~~10~~ & ~~0.8992E-03 & ~~0.9301E-03& ~~0.9332E-03      \\ 
    $L_2$ & ~~20 ~~& ~~0.1311E-04 & ~~0.3159E-04& ~~0.7893E-03    \\ 
    $$ & ~~50~~& ~~0.5141E-09 & ~~0.1361E-05& ~~0.3116E-03  \\

		\hline
		\\
		$$ &~~10~~ & ~~0.5275E-03 & ~~0.5497E-03& ~~0.5402E-03      \\ 
    $L_\infty$ & ~~20 ~~& ~~0.5221E-05 & ~~0.8033E-04& ~~0.3102E-04    \\ 
    $$ & ~~50~~& ~~0.2580E-09 & ~~0.3222E-06& ~~0.7445E-06  \\ 
		
		\hline
		\end{tabular}

\end{center}
Table 5 shows the numerical results for increasing number of grid points. Because of the asymptotic behaviour of the problem for all grid points results have approximately accurateness is approximately same. 

\section{Conclusion }
In this study it is examined an efficient method to yield solutions of higher order singularly perturbation  problems numerically. The Chebyshev based differential quadrature method applied to perturbation problems to obtain approximated results. Here, we chose linear and nonlinear test problems to show the accurateness of the presented method. The results are compared with exact solutions, as long as we have, or asymptotic solutions to support the idea. It is realized that the results are well-matched when compared to the known solutions.
\\
As a future work, both linear and nonlinear problems, ode and pde, can be solved using Chebyshev polynomials. Also, one can consider system of singularly perturbed boundary value problems.

\newpage

\newpage
\end{document}